\documentclass[11pt,a4paper]{article}
 \usepackage{indentfirst, latexsym,bm}
 \usepackage{amsmath,amssymb}
 \usepackage{mdframed}
 \usepackage{pifont}
 \usepackage{amsfonts}
 \usepackage{mathrsfs}
 \usepackage{array}
 \usepackage{multirow} 
 \usepackage{graphicx}
 \usepackage{subfigure}
 \usepackage{picinpar}
 \usepackage{setspace}
 \usepackage[textsize=footnotesize]{todonotes}
 \RequirePackage[colorlinks,citecolor=blue,urlcolor=blue,linkcolor=blue]{hyperref}
 \usepackage[top=2cm,bottom=2cm, outer=2cm, inner=2cm]{geometry}
 
 \usepackage{caption}
 \usepackage{subcaption}
 \usepackage[labelfont={bf,small},textfont={small}]{caption}
 
 \usepackage{tikz}
 \usetikzlibrary{arrows.meta}
 
 \RequirePackage[numbers]{natbib}
 \usepackage{enumitem} 
 
 \usepackage{booktabs}
 \usepackage{authblk}
 
 \usepackage{threeparttable}
 \usepackage{mathrsfs,amsfonts,amsmath}
 
 \usepackage{color}
 
 \makeatletter
 \def\namedlabel#1#2{\begingroup
 	#2%
 	\def\@currentlabel{#2}%
 	\phantomsection\label{#1}
 	\endgroup}
 \makeatother

 \numberwithin{figure}{section}

 \newcommand\email[1]{\href{mailto:#1}{ \nolinkurl{#1}}}

 \newtheorem{theorem}{Theorem}[section]
 \newtheorem{definition}[theorem]{Definition}
 \newtheorem{lemma}[theorem]{Lemma}
 \newtheorem{corollary}[theorem]{Corollary}
 \newtheorem{proposition}[theorem]{Proposition}
 \newtheorem{remark}[theorem]{Remark}
 \newtheorem{condition}[theorem]{Condition}
 \newtheorem{example}{Example}[section]

 \def\blemma{\begin{lemma}}\def\elemma{\end{lemma}}
 \def\bproposition{\begin{proposition}}\def\eproposition{\end{proposition}}
 \def\ttheorem{\begin{theorem}}\def\etheorem{\end{theorem}}
 \def\bcorollary{\begin{corollary}}\def\ecorollary{\end{corollary}}
 \def\bremark{\begin{remark}}\def\eremark{\end{remark}}
 \def\bcondition{\begin{condition}}\def\econdition{\end{condition}}

 \def\benumerate{\begin{enumerate}}\def\eenumerate{\end{enumerate}}
 \def\bitemize{\begin{itemize}}\def\eitemize{\end{itemize}}

 \def\beqlb{\begin{eqnarray}}\def\eeqlb{\end{eqnarray}}
 \def\beqnn{\begin{eqnarray*}}\def\eeqnn{\end{eqnarray*}}
 \def\ar{\!\!\!&}

 \def\qed{\hfill$\Box$\medskip}

\begin{document} 
 	\title{\bf  Compensation Formulas for Branching L\'evy Processes}
 	
 	\author{Hao Wu\footnote{School of Mathematical Sciences, Nankai University, China; email: wuhao.math@outlook.com} 
 		\quad\ and\quad 
 		Wei Xu\footnote{School of Mathematics and Statistics, Beijing Institute of Technology,  China; email: xuwei.math@gmail.com.}  }  
 	\maketitle

 \begin{abstract}
 This paper aims to establish two compensation formulas for branching L\'evy processes.  
 The first expresses expected sums over spatial jumps as integrals against the L\'evy measure, with the future family of the jumping particle started from its position immediately after the jump. 
 The second expresses expected sums over reproduction events as time integrals weighted by the branching rate and the offspring distribution, with each event marked by the ordered tuple of branching processes initiated by its offspring.  
 Combined with the fixed-time many-to-one formula, it recovers the whole-tree many-to-one identity for functionals of ancestral spatial paths.
 	
 	
 \end{abstract}

 \section{Introduction} 
 \label{Sec.Introduction}
 \setcounter{equation}{0}
 
 Consider a one-dimensional  branching L\'evy process with local reproduction defined on a probability space $\big(\Omega,\mathscr{F},\mathbf{P}\big)$. 
 The system starts with a single particle that moves according to a L\'evy process with triplet $(\beta,\sigma^2,\nu)$. 
 Its lifetime is exponentially distributed with rate $\gamma>0$ and is independent of its motion. 
 At the end of its lifetime, the particle is replaced by $k$ offspring with probability $p_k$, where $\boldsymbol{p}=\{p_k\}_{k\geq 0}$ is the offspring distribution. 
 All offspring are born at the parent's position at death and evolve independently according to the same motion and reproduction rules.
 We assume that only finitely many branching events occur in each bounded
 time interval almost surely. 
 
 This branching particle system  can also be constructed by equipping a continuous-time Galton--Watson tree with independent L\'evy motions along the lifetimes of its individuals; see \cite{BansayeDelmasMarsalleTran2011} and Section~\ref{Sec.Proofs}. 
 We label the initial particle by $\varnothing$ and other particles using the Ulam--Harris system. 
 The collection  $\mathbf{T} $ of all labels forms a subset of $ \mathcal{U}$ defined by
 \beqnn
 \mathcal{U} := \bigcup_{n=0}^{\infty}\mathbb{N}_+^n,
 \eeqnn
 with  $\mathbb{N}_+=\{1,2,\ldots\}$ and  the convention $\mathbb{N}_+^0=\{\varnothing\}$; see \cite{BansayeDelmasMarsalleTran2011,LeGall2005}  for details.  
 Let $N_f(\mathbb{R})$ denote the space of finite counting measures on
 $\mathbb{R}$, equipped with the topology of weak convergence. Let
 $D_f$ be the space of
 $N_f(\mathbb{R})$-valued c\`adl\`ag paths, equipped with the Skorokhod
 $J_1$-topology. For $z\in\mathbb{R}$, write $\delta_z$ for the unit
 point mass at $z$.
 For $t\geq 0$, denote by $V_t$ the set of
 particles alive at time $t$. 
 For $u\in V_t$ and $0\leq s\leq t$, let $X_u(s)$ be the position of $u$ or its ancestor at time $s$.  
 The $N_f(\mathbb{R})$-valued process $Z=\{Z_t \}_{t\geq 0} $ defined by
 \beqlb\label{eqn.101}
 Z_t:= \sum_{u\in V_t} \delta_{X_u(t)}
 \eeqlb
 has sample paths in $D_f$ and records the positions of all living particles over time. It is usually referred to as a branching L\'evy process. 
 We denote its natural filtration after the usual completion by $ \{\mathscr{F}_t\}_{t\geq 0}$. 
 For $x\in \mathbb{R}$, we write $\mathbf{P}_x$ and $\mathbf{E}_x$ for the law and expectation of this process with the initial particle started from $x\in\mathbb{R}$, and abbreviate $\mathbf{P}=\mathbf{P}_0$ and
 $\mathbf{E}=\mathbf{E}_0$.
 
 Our first result extends the L\'evy process compensation formula to the
 spatial jumps of all particles in the branching L\'evy process. 
 For $t\geq 0$ and $u \in V_t$, let $\Delta X_u(s):=X_u(s)-X_u(s-)$ denotes the jump of $X_u$ at time $s\in(0,t]$, and $Z^{t,u}=\{Z^{t,u}_s\}_{s\geq 0}$ denotes the future family of $u$, with its time origin at $t$ and its positions in the original spatial coordinates, i.e., 
 \beqnn
 Z^{t,u}_0=\delta_{X_u(t)}
 \quad \mbox{and}\quad 
 Z^{t,u}_s=\sum_{\substack{v\in V_{t+s}, u\preceq v}}
 \delta_{X_v(t+s)},\quad s> 0,
 \eeqnn 
  where $u\preceq v$ means that $u=v$ or $u$ is an ancestor of $v$. 
  In particular, $Z^{g_u,u}$ is the sub branching L\'evy process initiated by
 $u$ at its birth. 
 
 \begin{theorem}\label{MainThm01}
 Given a $(\mathscr{F}_t)$-predictable process $\boldsymbol{F}=\{ \boldsymbol{F}_t:t\geq 0 \}$ taking values in the space of non-negative measurable functions on  $\mathcal{U}\times\mathbb{R}$ such that $\boldsymbol{F}_t(u,0)=0$ for all $t\geq 0$ and $u\in \mathcal{U}$, and a deterministic measurable function $\boldsymbol{G}=\{ \boldsymbol{G}_t:t\geq 0 \}$ taking values in the space of non-negative measurable functions on  $\mathcal{U}\times \mathbb{R}\times  D_f$. For any  $x\in \mathbb{R}$,  we have 
 \beqlb\label{eqn.MainThm01}
  \lefteqn{\mathbf{E}_x\bigg[ \sum_{t>0} \sum_{\substack{u\in V_{t}}} \boldsymbol{F}_t\big(u,\Delta X_u(t)\big) \boldsymbol{G}_t\big(u, \Delta X_u(t), Z^{t,u}\big) \bigg] }\cr\cr
  \ar=\ar \int_0^\infty \mathbf{E}_x\bigg[   \sum_{u\in V_{t}}    \int_{\mathbb{R}} \boldsymbol{F}_t(u,y) \cdot \mathbf{E}_{X_u(t)+y}\big[ \boldsymbol{G}_t\big(u,y, Z\big) \big] \,\nu(dy) \bigg]\,dt . 
 \eeqlb
 \end{theorem}
 
 \begin{remark}
   When $p_1=1$, exactly one particle is alive at each time, and its position evolves as a L\'evy process. 
  When $p_0=1$, the initial particle dies without offspring, so its motion is killed at an independent exponential time of rate $\gamma$. 
  In these two cases, the formula \eqref{eqn.MainThm01} with $\boldsymbol{G}\equiv 1$ recovers the compensation formulas for a L\'evy process and for a L\'evy process killed at rate $\gamma$, respectively.
 	\end{remark}

 
 For $n\ge 1$, let $D_f^n$ be the space of ordered $n$-tuples of paths in $D_f$ endowed with the product Borel $\sigma$-algebra. 
 Set
 \beqnn
 \mathcal{D}_f:= \bigcup_{n=0}^{\infty}  D_f^n
 \quad\mbox{with}\quad
 \mathcal{D}_f^0:=\{()\},
 \eeqnn
 and equip $\mathcal{D}_f$ with the disjoint-union $\sigma$-algebra.  Here $()$ denotes the empty tuple. 
 For $x\in\mathbb{R}$ and $k\geq 1$, let $\mathbf{P}_x^{\otimes k}$
 denote the joint law on $D^k_f$ of $k$ independent copies $Z^1,\ldots,Z^k$ of
 the branching L\'evy process started from a single particle at $x$.
 Write $\mathbf{E}_x^{\otimes k}$ for the corresponding expectation.
 For $k=0$, let $\mathbf{P}_x^{\otimes 0}$ be the unit mass at $()$.

 Our second result gives a compensation formula for the reproduction events
 of the branching L\'evy process. 
 For each individual $u\in\mathbf{T}$, let $g_u$, $d_u$,
 and $\ell_u$ denote its birth time, death time, and lifetime, respectively. 
 We make the convention that $g_\varnothing=0$. 
 At its death time $d_u$, particle $u$ produces $k_u$ offspring, where
 $k_u$ has law $\boldsymbol{p}$. 
 Under the Ulam--Harris labeling, the $i$-th child of $u$ is labeled $ui$ for $1\leq i\leq k_u$.

 
 \begin{theorem}\label{MainThm02}
 Given a $(\mathscr{F}_t)$-predictable process $\boldsymbol{F}=\{ \boldsymbol{F}_t:t\geq 0 \}$ taking values in the space of non-negative measurable functions on $\mathcal{U}\times\mathbb{N}$, and  a  deterministic measurable function $\boldsymbol{G}=\{ \boldsymbol{G}_t :t\geq 0 \}$ taking values in the space of non-negative measurable functions on  $\mathcal{U}\times \mathcal{D}_f$.  For any  $x\in \mathbb{R}$,  we have
 \beqlb\label{eqn.MainThm02}
  \lefteqn{\mathbf{E}_x\bigg[	\sum_{u\in \mathbf{T}}  \boldsymbol{F}_{d_u}(u,k_u)\,  \boldsymbol{G}_{d_u}\Big(u,\big(Z^{d_u,ui}\big)_{i=1}^{k_u}\Big) \bigg]}\cr\cr
  \ar=\ar \mathbf{E}_x\bigg[	\sum_{t>0}\sum_{ u\in V_{t-}} \boldsymbol{F}_t(u,k_u)\cdot \mathbf{1}_{\{d_u=t\}}\cdot  \boldsymbol{G}_t\Big(u,\big(Z^{t,ui}\big)_{i=1}^{k_u}\Big) \bigg] \cr
  \ar=\ar \gamma\int_0^\infty  \mathbf{E}_x\bigg[
 		\sum_{u\in V_{t}}\sum_{k=0}^{\infty}p_k\,\boldsymbol{F}_t(u,k)\,
 		\mathbf{E}_{X_u(t)}^{\otimes k}
 		\Big[\boldsymbol{G}_t\Big(u,\big(Z^i\big)_{i=1}^k\Big)\Big] \bigg] dt. 
 \eeqlb
 \end{theorem}
 
 \begin{remark}
 The proof of Theorem~\ref{MainThm02} remains valid when the L\'evy motion
 is replaced by a strong Markov motion with a Borel measurable
 path law as a function of the starting point.	
 In particular, let $f$ be a non-negative measurable function on
 $\mathbb{R}_+\times D([0,\infty),\mathbb{R})$ such that $f(t,\omega)$
 depends only on $\omega|_{[0,t)}$. We write $f(t,\omega|_{[0,t)})$
 for its value. Taking $\boldsymbol{G}\equiv 1$, $ \boldsymbol{F}_t(u,k)
 =f\bigl(t,X_u|_{[0,t)}\bigr)$ for $u\in V_{t-}, k\geq 0$, and setting $\boldsymbol{F}_t(u,k)=0$ elsewhere, we obtain from
 \eqref{eqn.MainThm02} that
 \beqnn
 \mathbf{E}_x\bigg[\sum_{u\in\mathbf{T}}
 f\bigl(d_u,X_u|_{[0,d_u)}\bigr)\bigg]
 \ar=\ar
 \gamma\int_0^\infty \mathbf{E}_x\bigg[
 \sum_{u\in V_{t}}f\bigl(t,X_u|_{[0,t)}\bigr)\bigg]\,dt.
 \eeqnn
 If $m:=\sum_{k\geq 0}kp_k<\infty$, applying the fixed-time many-to-one
 formula to the right-hand side gives
 \beqnn
 \mathbf{E}_x\bigg[\sum_{u\in\mathbf{T}}
 f\bigl(d_u,X_u|_{[0,d_u)}\bigr)\bigg]
 \ar=\ar
 \gamma\int_0^\infty e^{\gamma(m-1)t} \cdot 
 \mathbf{E}_x^\xi\bigl[f(t,\xi|_{[0,t)})\bigr]\,dt,
 \eeqnn
 where $\mathbf{E}_x^\xi$ denotes expectation for a single copy $\xi$
 of the spatial motion started at $x$. This is the whole-tree many-to-one
 formula of \cite[Proposition~3.5]{BansayeDelmasMarsalleTran2011}
 specialized to local reproduction and functionals of the spatial path.
 	\end{remark}

 \begin{remark}
 	For any $x_0\in \mathbb{R}$, consider a branching L\'evy process $\widetilde{Z}$ with particles being killed as they enter $(-\infty,x_0)$. The two compensation formulas \eqref{eqn.MainThm01} and \eqref{eqn.MainThm02} remain valid for $\widetilde{Z}$, since $\widetilde{Z}$ can be constructed as a functional of $Z$, i.e., for any $t$ and $u\in \widetilde{V}_t$, we have $\widetilde{X}_u(t)=X_u(t)$ and 
 	\beqnn
 	 \widetilde{Z}_t =  \sum_{u\in V_{t}} \mathbf{1}_{\{ \inf\{ X_u(s):s\in[0,t] \} \geq x_0  \}}\cdot \delta_{X_u(t)}
 	 \quad \mbox{and}\quad 
 	  \widetilde{Z}^{t,u}_s=\sum_{\substack{v\in  \widetilde{V}_{t+s}, u\preceq v}}
 	\mathbf{1}_{\{\inf_{r\in[t,t+s]}X_v(r)\geq x_0  \}}\cdot \delta_{X_v(t+s)}, \quad s\geq0 . 
 	\eeqnn
 	Let ${\widetilde{\mathscr F}t}{t\geq0}$ denote the natural filtration of $\widetilde Z$, completed in the usual way.
 	Obviously, it is a subfiltration of $\{\mathscr{F}_t\}_{t\geq0}$. 
 	As an illustration, let $\boldsymbol F$ and $\boldsymbol G$ be two functionals satisfying the assumptions of Theorem~\ref{MainThm01}, with $\mathscr F_t$ replaced by $\widetilde{\mathscr F}_t$. 
 	Then 
 \beqnn
 		\lefteqn{
 			\mathbf{E}_x\Big[
 			\sum_{t>0}\sum_{u\in\widetilde{V}_t}
 			\boldsymbol{F}_t\big(u,\Delta\widetilde{X}_u(t)\big)
 			\boldsymbol{G}_t\big(u,\Delta\widetilde{X}_u(t),
 			\widetilde{Z}^{t,u}\big)
 			\Big]
 		}\cr\cr
 		\ar=\ar
 		\mathbf{E}_x\Big[
 		\sum_{t>0}\sum_{u\in V_t}
 		\mathbf{1}_{\{\underline{X}_u(t)\geq x_0\}}
 		\boldsymbol{F}_t\big(u,\Delta X_u(t)\big)
 		\boldsymbol{G}_t\big(u,\Delta X_u(t),
 		\widetilde{Z}^{t,u}\big)
 		\Big]\cr\cr
 		\ar=\ar
 		\int_0^\infty\mathbf{E}_x\Big[
 		\sum_{u\in V_{t-}}\int_{\mathbb{R}}
 		\mathbf{1}_{\{\underline{X}_u(t-)\geq x_0\}}
 		\mathbf{1}_{\{X_u(t-)+y\geq x_0\}}
 		\boldsymbol{F}_t(u,y)
 		\mathbf{E}_{X_u(t-)+y}
 		\big[\boldsymbol{G}_t(u,y,\widetilde{Z})\big]
 		\,\nu(dy)
 		\Big]\,dt\cr\cr
 		\ar=\ar
 		\int_0^\infty\mathbf{E}_x\Big[
 		\sum_{u\in\widetilde{V}_t}\int_{\mathbb{R}}
 		\mathbf{1}_{\{\widetilde{X}_u(t)+y\geq x_0\}}
 		\boldsymbol{F}_t(u,y)
 		\mathbf{E}_{\widetilde{X}_u(t)+y}
 		\big[\boldsymbol{G}_t(u,y,\widetilde{Z})\big]
 		\,\nu(dy)
 		\Big]\,dt.
 		\eeqnn
 	\end{remark}

 \section{Proofs} 
 \label{Sec.Proofs}
 \setcounter{equation}{0}

 Before proving the two main theorems, we first describe in detail the construction of the branching L\'evy process $Z$ started from a single particle at $x\in\mathbb{R}$ by assigning independent L\'evy motions to the individuals of a continuous-time Galton--Watson tree.
 To each potential label $u\in\mathcal{U}$, we assign an exponential lifetime $\ell_u$ of rate $\gamma$, an offspring count $k_u$ with law $\boldsymbol{p}$, and a L\'evy process $Y_u=(Y_u(s))_{s\geq 0}$ started at zero with triplet $(\beta,\sigma^2,\nu)$.
 The parameter $s$ denotes the age of the particle.
 All these random elements are mutually independent, both within and across labels. 
 
 The realized tree $\mathbf{T}$ contains the root $\varnothing$. For a
 label $u=(u_1,\ldots,u_n)$ with $n\geq 1$,  it belongs to $\mathbf{T}$ if and only if $u_i\leq k_{u^{(i-1)}}$ for all $i=1,\ldots,n$, where  $u^{(0)}=\varnothing$ and $u^{(i-1)}=(u_1,\ldots,u_{i-1})$ for $i\geq 2$.  
 For each $u\in \mathbf{T}$, its birth  and death times are
 \beqnn
 g_u:= \sum_{v \prec u} \ell_v
 \quad \mbox{and}\quad 
 d_u:= g_u+\ell_u, 
 \eeqnn
 where $v\prec u$ means that $v$ is a strict ancestor of $u$. 
 The empty sum gives $g_\varnothing=0$. 
 Thus the set of particles alive at time $t\geq 0$ is  $V_t= \{ u\in \mathbf{T}: g_u\leq t<d_u \}$. 
 For $u\in V_t$, its position at time $t$ can be given by 
  \beqnn
X_u(t)= x+ \sum_{v \prec u} Y_v(\ell_v) + Y_u(t-g_u). 
\eeqnn 
 Each strict ancestor contributes its displacement over its entire lifetime, and $u$ contributes its displacement up to age $t-g_u$.
 Substituting $V_t$ and $X_u(t)$ into \eqref{eqn.101} gives a realization of $Z$.
 
 \medskip
 
   \textit{\textbf{Proof of Theorem~\ref{MainThm01}.}} 
  Let $\mathcal{T}=\{(u,\ell_u) \}_{u\in \mathbf{T}}$ be the whole genealogy of all labels and lifetimes.  
  Almost surely, no jump occurs at $g_u$ for all $u\in\mathbf{T}$, since the lifetime is independent of the motion and a L\'evy process has no jump at a fixed age.  
  Reindexing all jumps gives
  \beqlb \label{eqn.301}
  S\ar:=\ar\sum_{t>0} \sum_{\substack{u\in V_t}} \boldsymbol{F}_t\big(u,\Delta X_u(t)\big) \, \boldsymbol{G}_t\big(u,\Delta X_u(t), Z^{t,u}\big) \cr\cr
  \ar=\ar \sum_{u\in \mathbf{T}}  \sum_{t\in [g_u,d_u)} \boldsymbol{F}_t\big(u,\Delta X_u(t)\big)\, \boldsymbol{G}_t\big(u,\Delta X_u(t), Z^{t,u}\big). 
  \eeqlb 
  Taking expectations on  both sides of the second equality, then using the
  tower property and conditional Tonelli theorem, we obtain
  \beqlb\label{eqn.303}
  \mathbf{E}_x[S]  
  \ar=\ar  \mathbf{E}_x\Big[  \sum_{u\in \mathbf{T}} S^\mathcal{T}_u \Big] 
  \quad \mbox{with}\quad S^\mathcal{T}_u:=\mathbf{E}_x\bigg[   \sum_{t\in [g_u,d_u)} \boldsymbol{F}_t\big(u,\Delta X_u(t)\big)\, \boldsymbol{G}_t\big(u,\Delta X_u(t),Z^{t,u}\big)\,\bigg|\, \mathcal{T}\bigg]. 
  \eeqlb
  For $\epsilon>0$ and $u\in\mathbf{T}$, set $\tau_{u,0}^{\epsilon}=g_u$ and define recursively $\tau_{u,j}^{\epsilon}:= \inf\{ s\in(\tau_{u, j-1}^{\epsilon},d_u): |\Delta X_u(s)| \geq \epsilon \}$ for $j\geq 1$ with the convention $\inf\emptyset =\infty$. 
  The finite terms of $\{\tau_{u,j}^{\epsilon}\}_{j\geq 1}$ are the successive
  times at which particle $u$ makes an upward or downward jump of size at least
  $\epsilon$ during its lifetime. 
  Let
  \beqnn
  S^\mathcal{T}_u(\epsilon)
  \ar:=\ar \mathbf{E}_x\bigg[  \sum_{j=1}^\infty \mathbf{F}_{\tau_{u,j}^{\epsilon}}\big(u,\Delta X_u(\tau_{u,j}^{\epsilon})\big)\, \mathbf{G}_{\tau_{u,j}^{\epsilon}}\big(u,\Delta X_u(\tau_{u,j}^{\epsilon}),Z^{\tau_{u,j}^{\epsilon},u}\big)  \,\bigg|\, \mathcal{T}\bigg]. 
  \eeqnn
  The monotone convergence theorem yields that $S^\mathcal{T}_u(\epsilon) \overset{\rm a.s.}\to S^\mathcal{T}_u$ as $\epsilon\to 0+$. 
  For $t\geq 0$, let $\mathscr{G}_t:= \mathscr{F}_t\vee  \sigma\big( \mathcal{T} \big)$ reveal $\mathcal{T}$ at time zero and each particle's
  spatial motion up to time $t$.   
  It is obvious that $\{\tau_{u,j}^\epsilon\}_{j\geq 1}$ is a sequence of  $(\mathscr{G}_t)$-stopping times.  
  Since  $\sigma(\mathcal{T})\subset\mathscr{G}_{\tau^\epsilon_{u,1}}$, we have 
  \beqlb\label{eqn.302}
  S_u^\mathcal{T}(\epsilon) 
  \ar=\ar \mathbf{E}_x\big[ \mathbf{E}_x[ S_u^\mathcal{T}(\epsilon)\,|\, \mathscr{G}_{\tau^\epsilon_{u,1}}] \,\big|\,\mathcal{T} \big] \cr\cr
  \ar=\ar \mathbf{E}_x\bigg[ \mathbf{E}_x\Big[ \sum_{j=1}^\infty \mathbf{F}_{\tau_{u,j}^{\epsilon}}\big(u,\Delta X_u(\tau_{u,j}^{\epsilon})\big)\, \mathbf{G}_{\tau_{u,j}^{\epsilon}}\big(u,\Delta X_u(\tau_{u,j}^{\epsilon}),Z^{\tau_{u,j}^{\epsilon},u}\big) \,\Big|\, \mathscr{G}_{\tau^\epsilon_{u,1}} \Big] \bigg| \mathcal{T}\bigg] . 
  \eeqlb  
  Conditioned on a given genealogy $\mathtt{T}$, for any particle $u$ and  $t \in [g_u,d_u)$, let $Q^\mathtt{T}_{t,u}(x,dz)$ be the law on $D_f$ of  branching L\'evy process started at $x \in \mathbb{R}$ at time $t$, using the branching schedule of $\mathtt{T}$ from time $t$ and fresh independent L\'evy motions. 
  Define
  \beqlb\label{eqn.304}
  \boldsymbol{K}^{\mathtt{T}}_{t,u}(x,y)
  :=\int_{\mathcal  D_f}\boldsymbol{G}_t(u,y,z)\,Q^{\mathtt{T}}_{t,u}(x,dz),\quad x,y\in \mathbb{R}.
  \eeqlb
  On the event $\tau_{u,1}^\epsilon \in [g_u,d_u)$,   we have $\boldsymbol{F}_{\tau_{u,1}^{\epsilon}}\big(u,\Delta X_u(\tau_{u,1}^{\epsilon})\big)$ is $\mathscr{G}_{\tau^\epsilon_{u,1}}$-measurable and
  \beqnn
  \mathbf{E}_x\Big[  \boldsymbol{G}_{\tau_{u,1}^{\epsilon}}\big(u,\Delta X_u(\tau_{u,1}^{\epsilon}), Z^{\tau_{u,1}^{\epsilon},u}\big) \,\Big|\, \mathscr{G}_{\tau^\epsilon_{u,1}} \Big] 
  \ar=\ar   \boldsymbol{K}^{\mathcal{T}}_{\tau_{u,1}^{\epsilon},u}\big(X_u(\tau_{u,1}^{\epsilon}),\Delta X_u(\tau_{u,1}^{\epsilon})\big) .
  \eeqnn
  Substituting these into the last expression in \eqref{eqn.302} yields  
  \beqnn
  S^\mathcal{T}_u(\epsilon) 
  \ar=\ar  \mathbf{E}_x\bigg[   \boldsymbol{F}_{\tau_{u,1}^{\epsilon}}\big(u,\Delta X_u(\tau_{u,1}^{\epsilon})\big)\,  \boldsymbol{K}^{\mathcal{T}}_{\tau_{u,1}^{\epsilon},u}\big(X_u(\tau_{u,1}^{\epsilon}),\Delta X_u(\tau_{u,1}^{\epsilon})\big)  \, \bigg|\, \mathcal{T}\bigg]\cr
  \ar\ar  + \mathbf{E}_x\bigg[   \sum_{j=2}^\infty \mathbf{F}_{\tau_{u,j}^{\epsilon}}\big(u,\Delta X_u(\tau_{u,j}^{\epsilon})\big)\, \mathbf{G}_{\tau_{u,j}^{\epsilon}}\big(u,\Delta X_u(\tau_{u,j}^{\epsilon}), Z^{\tau_{u,j}^{\epsilon},u}\big)  \, \bigg|\, \mathcal{T}\bigg] . 
  \eeqnn
  Using the preceding argument recursively to the last conditional expectation gives
  \beqnn
  S^\mathcal{T}_u(\epsilon) 
  \ar=\ar  \mathbf{E}_x\bigg[   \sum_{j=1}^\infty \boldsymbol{F}_{\tau_{u,j}^{\epsilon}}\big(u,\Delta X_u(\tau_{u,j}^{\epsilon})\big)\,  \boldsymbol{K}^{\mathcal{T}}_{\tau_{u,j}^{\epsilon},u}\big(X_u(\tau_{u,j}^{\epsilon}),\Delta X_u(\tau_{u,j}^{\epsilon})\big)  \, \bigg|\, \mathcal{T}\bigg]\cr
  \ar=\ar  \mathbf{E}_x\bigg[   \sum_{t\in [g_u,d_u)}\mathbf{1}_{\{ |\Delta X_u(t)|\geq \epsilon \}}\cdot  \boldsymbol{F}_{t}\big(u,\Delta X_u(t)\big)\, \boldsymbol{K}^{\mathcal{T}}_{t,u}\big(X_u(t-)+\Delta X_u(t),\Delta X_u(t) \big)  \, \bigg|\, \mathcal{T}\bigg], 
  \eeqnn
  which  converges almost surely as $\epsilon\to 0+$ to 
  \beqnn
  S^\mathcal{T}_u 
  \ar=\ar  \mathbf{E}_x\bigg[   \sum_{t\in [g_u,d_u)}  \boldsymbol{F}_{t}\big(u,\Delta X_u(t)\big)\, \boldsymbol{K}^{\mathcal{T}}_{t,u}\big(X_u(t-)+\Delta X_u(t),\Delta X_u(t)\big)  \, \bigg|\, \mathcal{T}\bigg].
  \eeqnn
  Conditioned on $\mathcal{T}$, the process $X_u$ is a L\'evy process on $[0,d_u]$ started at $x$ with L\'evy measure $\nu$. 
  By the compensation formula for jumps of L\'evy process; see e.g. \cite[p.7]{Bertoin1996}, 
  \beqnn
  S^\mathcal{T}_u 
  \ar=\ar  \mathbf{E}_x\bigg[   \int_{g_u}^{d_u} dt \int_{\mathbb{R}} \boldsymbol{F}_{t}\big(u,y\big)\, \boldsymbol{K}^{\mathcal{T}}_{t,u}\big(X_u(t-)+y,y\big) \,\nu(dy) \, \bigg|\, \mathcal{T}\bigg].
  \eeqnn
  Substituting this into \eqref{eqn.303} gives
  \beqnn
  \mathbf{E}_x[S]  
  \ar=\ar  \mathbf{E}_x\bigg[  \sum_{u\in \mathbf{T}} \mathbf{E}_x\Big[   \int_{g_u}^{d_u} dt \int_{\mathbb{R}} \boldsymbol{F}_{t}\big(u,y\big)\, \boldsymbol{K}^{\mathcal{T}}_{t,u}\big(X_u(t-)+y,y\big) \,\nu(dy) \, \Big|\, \mathcal{T}\Big] \bigg] \cr
  \ar=\ar \mathbf{E}_x\bigg[  \sum_{u\in \mathbf{T}}  \int_{g_u}^{d_u} dt \int_{\mathbb{R}} \boldsymbol{F}_{t}\big(u,y\big)\, \boldsymbol{K}^{\mathcal{T}}_{t,u}\big(X_u(t-)+y,y\big) \,\nu(dy)   \bigg] \cr
  \ar=\ar \mathbf{E}_x\bigg[    \int_0^\infty dt \sum_{u\in V_t} \int_{\mathbb{R}} \boldsymbol{F}_{t}\big(u,y\big)\, \boldsymbol{K}^{\mathcal{T}}_{t,u}\big(X_u(t-)+y,y\big) \,\nu(dy)   \bigg] .
  \eeqnn
  Since $\boldsymbol{F}_t(u,y)$ is
  $\mathscr F_t$-measurable, by Tonelli's theorem and the tower property we have 
  \beqlb\label{eqn.305}
  \mathbf{E}_x[S]  \ar=\ar  \int_0^\infty dt \int_{\mathbb{R}}\mathbf{E}_x\bigg[   \sum_{u\in V_t}  \boldsymbol{F}_{t}\big(u,y\big)\, \boldsymbol{K}^{\mathcal{T}}_{t,u}\big(X_u(t-)+y,y\big)  \bigg] \,\nu(dy)  \cr
  \ar=\ar  \int_0^\infty dt \int_{\mathbb{R}}\mathbf{E}_x\bigg[   \sum_{u\in V_t}  \boldsymbol{F}_{t}\big(u,y\big)\, \mathbf{E}_x\Big[\boldsymbol{K}^{\mathcal{T}}_{t,u}\big(X_u(t-)+y,y\big) \,\Big|\, \mathscr{F}_t\Big] \bigg] \,\nu(dy)  . 
  \eeqlb
  By  the exponential memorylessness and the independent
  future motions and reproduction, the branching property at time $t$ yields that conditional on $\mathscr F_t$, the process $Z^{t,u}$ is a branching L\'evy process started at $X_u(t-)+y$. 
  This along with \eqref{eqn.304} induces that for $u\in V_t$, 
  \beqnn
 	\mathbf E_x\Big[
 	\boldsymbol K^{\mathcal T}_{t,u}\big(X_u(t-)+y,y\big)
 	\,\Big|\,\mathscr F_t\Big]
 	\ar=\ar
 	\mathbf E_x\Big[
 	\int_{D_f}
 	\boldsymbol G_t(u,y,\omega)\,
 	Q^{\mathcal T}_{t,u}\big(X_u(t-)+y,d\omega\big)
 	\,\Big|\,\mathscr F_t\Big]
 	\cr\cr
 	\ar=\ar
 	\int_{D_f}\boldsymbol G_t(u,y,\omega)\,
 	\mathbf E_x\Big[
 	Q^{\mathcal T}_{t,u}\big(X_u(t-)+y,d\omega\big)
 	\,\Big|\,\mathscr F_t\Big]
 	\cr\cr
 	\ar=\ar
 	\int_{D_f}\boldsymbol G_t(u,y,\omega)\,
 	\mathbf P_{X_u(t-)+y}(Z\in d\omega)=
 	\mathbf E_{X_u(t-)+y}
 	\big[\boldsymbol G_t(u,y,Z)\big].
 	\eeqnn
 The desired formula \eqref{eqn.MainThm01} follows directly by plugging this into the right-hand side of the second equality in \eqref{eqn.305} and using the stochastic continuity of L\'evy processes.
 This completes the proof. 
 \qed 
 
 \medskip
 
 \textit{\textbf{Proof of Theorem~\ref{MainThm02}.}} 
 The first equality in \eqref{eqn.MainThm02} follows directly from the fact that  each particle contributes to the sum only at its death time. 
 For the second equality, we first have  
 \beqlb\label{eqn.401}
 \lefteqn{\mathbf{E}_x\bigg[	\sum_{u\in \mathbf{T}}  \boldsymbol{F}_{d_u}(u,k_u)\,  \boldsymbol{G}_{d_u}\Big(u,\big(Z^{d_u,ui}\big)_{i=1}^{k_u}\Big) \bigg]}\cr\cr
 \ar=\ar \mathbf{E}_x\bigg[	  \boldsymbol{F}_{d_\varnothing}(\varnothing,k_\varnothing)\,  \boldsymbol{G}_{d_\varnothing}\Big(\varnothing,\big(Z^{d_\varnothing,i}\big)_{i=1}^{k_\varnothing}\Big) \bigg]
  + \mathbf{E}_x\bigg[	\sum_{u\in \mathbf{T}, u\neq \varnothing}  \boldsymbol{F}_{d_u}(u,k_u)\,  \boldsymbol{G}_{d_u}\Big(u,\big(Z^{d_u,ui}\big)_{i=1}^{k_u}\Big) \bigg] . 
 \eeqlb
 For the first expectation of the right-hand side, since $\boldsymbol{F}_{d_\varnothing}(\varnothing,k_\varnothing)$ is $\mathscr{F}_{d_\varnothing}$-measurable, by the tower property and the branching property we have  
 \beqnn
 \mathbf{E}_x\bigg[	  \boldsymbol{F}_{d_\varnothing}(\varnothing,k_\varnothing)\,  \boldsymbol{G}_{d_\varnothing}\Big(\varnothing,\big(Z^{d_\varnothing,i}\big)_{i=1}^{k_\varnothing}\Big) \bigg] 
 \ar=\ar 
 \mathbf{E}_x\bigg[	\boldsymbol{F}_{d_\varnothing}(\varnothing,k_\varnothing)\cdot  \mathbf{E}_x\Big[   \boldsymbol{G}_{d_\varnothing}\Big(\varnothing,\big(Z^{d_\varnothing,i}\big)_{i=1}^{k_\varnothing}\Big) \,\Big|\, \mathscr{F}_{d_\varnothing}\Big] \bigg].
 \eeqnn
The parent's motion has no jump at its independent exponential death
age almost surely. Hence every child of $\varnothing$ starts at $X_\varnothing(d_\varnothing-)$. 
Consequently, 
\beqnn
\mathbf{E}_x\Big[   \boldsymbol{G}_{d_\varnothing}\Big(\varnothing,\big(Z^{d_\varnothing,i}\big)_{i=1}^{k_\varnothing}\Big) \,\Big|\, \mathscr{F}_{d_\varnothing}\Big]
\ar=\ar \mathbf E_{X_\varnothing(d_\varnothing-)}^{\otimes k_\varnothing}\Big[   \boldsymbol{G}_{d_\varnothing}\Big(\varnothing,\big(Z^{i}\big)_{i=1}^{k_\varnothing}\Big)\Big],
\eeqnn
which  follows that
 \beqlb\label{eqn.403}
 \mathbf{E}_x\bigg[	  \boldsymbol{F}_{d_\varnothing}(\varnothing,k_\varnothing)\,  \boldsymbol{G}_{d_\varnothing}\Big(\varnothing,\big(Z^{d_\varnothing,i}\big)_{i=1}^{k_\varnothing}\Big) \bigg]
 \ar=\ar \mathbf{E}_x\bigg[	\boldsymbol{F}_{d_\varnothing}(\varnothing,k_\varnothing)\cdot  \mathbf E_{X_\varnothing(d_\varnothing-)}^{\otimes k_\varnothing}\Big[   \boldsymbol{G}_{d_\varnothing}\Big(\varnothing,\big(Z^{i}\big)_{i=1}^{k_\varnothing}\Big)\Big] \bigg]. 
 \eeqlb
 For $ i=1,\ldots, k_\varnothing$, let $\mathbf{T}^{i}$ collects labels of the $i$-th child of $\varnothing$ and all its  descendants. 
 By the tower property and the mutual independence among descendants, the second expectation on the right-hand side of \eqref{eqn.401},
 \beqlb\label{eqn.402}
  \lefteqn{\mathbf{E}_x\bigg[ \mathbf{E}_x\Big[\sum_{u\in \mathbf{T}, u\neq \varnothing}  \boldsymbol{F}_{d_u}(u,k_u)\,  \boldsymbol{G}_{d_u}\Big(u,\big(Z^{d_u,uj}\big)_{j=1}^{k_u}\Big)\,\Big|\,\mathscr{F}_{d_\varnothing} \Big]\bigg] }\cr\cr
  \ar=\ar \mathbf{E}_x\bigg[ \sum_{i=1}^{k_\varnothing} \mathbf{E}_x\Big[\sum_{u\in \mathbf{T}^{i}}  \boldsymbol{F}_{d_u}(u,k_u)\,  \boldsymbol{G}_{d_u}\Big(u,\big(Z^{d_u,uj}\big)_{j=1}^{k_u}\Big)\,\Big|\,\mathscr{F}_{d_\varnothing} \Big]\bigg] . 
 \eeqlb
 Similarly as in the preceding argument, for any $ i=1,\ldots, k_\varnothing$, we have 
 \beqnn
 \lefteqn{\mathbf{E}_x\Big[\sum_{u\in \mathbf{T}^{i}}  \boldsymbol{F}_{d_u}(u,k_u)\,  \boldsymbol{G}_{d_u}\Big(u,\big(Z^{d_u,ui}\big)_{i=1}^{k_u}\Big)\,\Big|\,\mathscr{F}_{d_\varnothing} \Big]}\cr\cr
 \ar=\ar \mathbf{E}_x\Big[  \boldsymbol{F}_{d_{i}}(i,k_i)\,  \boldsymbol{G}_{d_{ i}}\Big(i,\big(Z^{d_{i}, ij}\big)_{j=1}^{k_{i}}\Big)\,\Big|\,\mathscr{F}_{d_\varnothing} \Big]
 + \mathbf{E}_x\Big[\sum_{u\in \mathbf{T}^{i},u\neq i}  \boldsymbol{F}_{d_u}(u,k_u)\,  \boldsymbol{G}_{d_u}\Big(u,\big(Z^{d_u,uj}\big)_{j=1}^{k_u}\Big)\,\Big|\,\mathscr{F}_{d_\varnothing} \Big]
 \eeqnn
 and further the first expectation on the right-hand side equals
 \beqnn
  \mathbf{E}_x\bigg[  \boldsymbol{F}_{d_i}(i,k_i)\,  \mathbf E_{X_{i}(d_{i}-)}^{\otimes k_{i}}\Big[   \boldsymbol{G}_{d_{i} }\Big(i,\big(Z^{j}\big)_{ j=1}^{k_{i}}\Big)\Big]\,\bigg|\,\mathscr{F}_{d_\varnothing} \bigg]
 \eeqnn
 Substituting these back into the right-hand side of \eqref{eqn.402} gives 
 \beqnn
  \lefteqn{\mathbf{E}_x\bigg[ \mathbf{E}_x\Big[\sum_{u\in \mathbf{T}, u\neq \varnothing}  \boldsymbol{F}_{d_u}(u,k_u)\,  \boldsymbol{G}_{d_u}\Big(u,\big(Z^{d_u,ui}\big)_{i=1}^{k_u}\Big)\,\Big|\,\mathscr{F}_{d_\varnothing} \Big]\bigg] }\cr\cr
  \ar=\ar \mathbf{E}_x\bigg[ \sum_{i=1}^{k_\varnothing} \boldsymbol{F}_{d_i}(i,k_i)\,  \mathbf E_{X_{i}(d_{i}-)}^{\otimes k_{i}}\Big[   \boldsymbol{G}_{d_{i} }\Big(i,\big(Z^{j}\big)_{ j=1}^{k_{i}}\Big)\Big] \bigg] +\mathbf{E}_x\bigg[ \sum_{i=1}^{k_\varnothing} \sum_{u\in \mathbf{T}^{ i},u\neq i}  \boldsymbol{F}_{d_u}(u,k_u)\,  \boldsymbol{G}_{d_u}\Big(u,\big(Z^{d_u,ui}\big)_{i=1}^{k_u}\Big) \bigg]. 
 \eeqnn 
 Using the preceding argument recursively to the last  expectation, then plugging this and \eqref{eqn.403} back into the right-hand sides of \eqref{eqn.401}, we obtain 
 \beqnn
  \mathbf{E}_x\bigg[	\sum_{u\in \mathbf{T}}  \boldsymbol{F}_{d_u}(u,k_u)\,  \boldsymbol{G}_{d_u}\Big(u,\big(Z^{d_u,ui}\big)_{i=1}^{k_u}\Big) \bigg]
  \ar=\ar \mathbf{E}_x\bigg[ \sum_{u\in \mathbf{T} } \boldsymbol{F}_{d_u}(u,k_u)\,  \mathbf E_{X_{u}(d_{u}-)}^{\otimes k_{u}}\Big[   \boldsymbol{G}_{d_{u} }\Big(u,\big(Z^{j}\big)_{ j=1}^{k_{u}}\Big)\Big] \bigg] . 
 \eeqnn
 
 We now continue considering the expectation on the right-hand side of the last equation. 
 For convenience, we set   
 \beqnn
 \boldsymbol{H}_t(u,k) := \boldsymbol{F}_{t}(u,k)\,  \mathbf E_{X_{u}(t-)}^{\otimes k}\Big[   \boldsymbol{G}_{t}\Big(u,\big(Z^{j}\big)_{ j=1}^{k}\Big)\Big],\quad t\geq 0,\, u\in V_{t-},\, k \geq 0, 
 \eeqnn
 and $\boldsymbol{H}_t(u,k)=0$ elsewhere. 
 The indicator of $\{u\in V_{t-}\}=\{g_u<t\leq d_u\}$ is $(\mathscr{F}_t)$-predictable, with $g_u=d_u=\infty$ for $u\notin\mathbf{T}$.  
 Thus $\boldsymbol{H}_t(u,k)$ is $(\mathscr{F}_t)$-predictable and  
 \beqnn
 \mathbf{E}_x\bigg[	\sum_{u\in \mathbf{T}}  \boldsymbol{F}_{d_u}(u,k_u)\,  \boldsymbol{G}_{d_u}\Big(u,\big(Z^{d_u,ui}\big)_{i=1}^{k_u}\Big) \bigg] 
 \ar=\ar \mathbf{E}_x\bigg[ \sum_{u\in \mathbf{T} }  \boldsymbol{H}_{d_u}(u,k_u)  \bigg] \cr\cr
 \ar=\ar \sum_{ u\in \mathcal{U}} \mathbf{E}_x\Big[  \mathbf{1}_{\{u\in \mathbf{T} \}} \cdot\boldsymbol{H}_{d_u}(u,k_u)  \Big].
 \eeqnn
 For each $u\in \mathcal{U}$, let $\mathscr{A}_u:= \sigma\big((Y_v)_{v\in\mathcal{U}},
 (\ell_v,k_v)_{v\in\mathcal{U}\setminus\{u\}}\big)$. 
 Hence both the event
 $\{u\in\mathbf{T}\}$ and the time $g_u$ are $\mathscr{A}_u$-measurable,
 since they depend only on the strict ancestors of $u$. 
 Conditionally on $\mathscr{A}_u$, the lifetime $\ell_u$ and the offspring count
 $k_u$ are still independent. 
 For $0\leq a<b<\infty$, $A\in\mathscr{F}_a$ and $k\geq 0$, we first consider the contribution of particle $u$ when $A$ occurs and $a<d_u\leq b$.
 On the event $\{u\in\mathbf{T},\,d_u>a\}$, particle $u$ has not died by time $a$,  so the history up to $a$ is the same if we replace $\ell_u$ by $a+1$
 and $k_u$ by zero while keeping all other variables fixed. 
  Consequently, on this event, $\mathbf{1}_A=\mathbf{1}_{A^u}$ for
 some $A^u\in\mathscr{A}_u$, obtained by evaluating $A$ in that
 modified history.
 Since $d_u=g_u+\ell_u$ on $\{u\in\mathbf{T}\}$, the density of $\ell_u$
 and the law of $k_u$ give
 \beqlb\label{eqn.404}
 \mathbf{E}_x\big[
 \mathbf{1}_A \cdot \mathbf{1}_{\{a<d_u\leq b,\,k_u=k\}}\,\big|\,\mathscr{A}_u\big]
 \ar=\ar \mathbf{1}_{A^u}\,p_k
 \int_{a\vee g_u}^{b}\gamma e^{-\gamma(t-g_u)}\,dt \cr\cr
 \ar=\ar \gamma p_k\,\mathbf{E}_x\bigg[
 \mathbf{1}_A\int_a^b \mathbf{1}_{\{g_u<t\leq d_u\}}\,dt
 \,\bigg|\,\mathscr{A}_u\bigg],
 \eeqlb
 where the middle integral is taken to be zero if $a\vee g_u\geq b$.
 Here the last equality follows from Tonelli's theorem and the fact that
 \beqnn
 \mathbf{P}_x \big( d_u \in dt ,k_u =k \,\big|\,\mathscr{A}_u \big) 
 =    p_k\cdot \gamma e^{-\gamma(t-g_u)} \cdot \mathbf{1}_{\{ g_u\leq t \}} \, dt,\quad t\geq 0, k\geq 0.
 \eeqnn 
 When $u\notin\mathbf{T}$, both conditional expectations in \eqref{eqn.404} are zero. 
 Taking expectations on both sides of \eqref{eqn.404} yields the desired identity for a predictable rectangle $(a,b]\times A$ with mark $k$. 
 For $A\in\mathscr{F}_a$, we apply the same identity on $(a+1/n,b]$ using an event in $\mathscr{F}_{a+1/n}$ that agrees with $A$ almost surely, and then
 let $n\to\infty$. 
 By the monotone class theorem and the monotone convergence theorem, we obtain
 \beqnn
 \mathbf{E}_x\Big[\mathbf{1}_{\{u\in\mathbf{T}\}}\,
 \boldsymbol{H}_{d_u}(u,k_u)\Big]
 \ar=\ar \gamma\,\mathbf{E}_x\bigg[
 \int_0^\infty\mathbf{1}_{\{g_u<t\leq d_u\}}
 \sum_{k=0}^\infty p_k\,\boldsymbol{H}_t(u,k)\,dt\bigg].
 \eeqnn
 Finally, summing this identity over $u\in\mathcal{U}$ and then using
 Tonelli's theorem, we have
 \beqnn
 \mathbf{E}_x\bigg[\sum_{u\in\mathbf{T}}\boldsymbol{H}_{d_u}(u,k_u)\bigg]
 \ar=\ar \gamma\int_0^\infty\mathbf{E}_x\bigg[
 \sum_{u\in V_{t-}}\sum_{k=0}^\infty p_k\,\boldsymbol{H}_t(u,k)
 \bigg]\,dt .
 \eeqnn
 The second equality in \eqref{eqn.MainThm02} follows directly by
 substituting the definition of $\boldsymbol{H}_t(u,k)$ into the
 right-hand side of the last equality, then using the stochastic continuity of $V$ and L\'evy processes. This completes the proof.
 \qed

  \bibliographystyle{plain}
 \bibliography{Reference}

@Article{BansayeDelmasMarsalleTran2011,
 	author  = {Bansaye, V. and Delmas, J.F. and  Marsalle, L. and Tran, V.C.},
 	title   = {{Limit theorems for Markov processes indexed by continuous time
 	Galton--Watson trees}},
 	year    = {2011},
 	journal = {Ann. Appl. Probab.},
 	volume  = {21},
 	number  = {6},
 	pages   = {2263--2314},
 }

@Book{Bertoin1996,
	author    = {Bertoin, J.},
	publisher = {Cambridge University Press},
	title     = {L\'evy {P}rocesses},
	year      = {1996},
	address   = {Cambridge},
}

@Article{LeGall2005,
 	author  = {Le Gall, J.F.},
 	title   = {Random trees and applications},
 	journal = {Probab. Surv.},
 	year    = {2005},
 	volume  = {2}, 
 	pages   = {245--311},
 }
 \end{document}